\documentclass{article}
\usepackage{amsmath}
\usepackage{amsthm}
\usepackage{amssymb}
\usepackage{amsfonts}
\usepackage{mathrsfs}
\usepackage{enumerate}
\usepackage{authblk}
\usepackage[normalem]{ulem} 
\usepackage[dvipsnames]{xcolor}
\usepackage[utf8]{inputenc}
\usepackage[english]{babel}
\usepackage[ top = 3 cm, bottom = 2 cm]{geometry} 
\usepackage{graphicx}

\usepackage[colorlinks,citecolor=red,pagebackref,hypertexnames=false]{hyperref}
\usepackage{pdfsync}

\numberwithin{equation}{section}
\theoremstyle{plain}
\newtheorem{theorem}{Theorem}
\newtheorem{lemma}[theorem]{Lemma}

\newtheorem{corollary}[theorem]{Corollary}

\theoremstyle{definition}
\newtheorem{definition}[theorem]{Definition}

\begin{document}

\title{Discretization and Convergence of the EIT Optimal Control Problem in 2D and 3D Domains}
\author{Ugur G. Abdulla\thanks{abdulla@fit.edu} \ } 
\author{Saleheh Seif\thanks{sseif2014@my.fit.edu}}
\affil{Department of Mathematical Sciences, Florida Institute of Technology,\\ Melbourne, Florida 32901}

\maketitle

\abstract{ We consider Inverse Electrical Impedance Tomography (EIT) problem on recovering electrical conductivity and potential in the body based on the measurement of the boundary voltages on the $m$ electrodes for a given electrode current. The variational formulation is pursued in the optimal control framework, where electrical conductivity and boundary voltages are control parameters, and the cost functional is the norm declinations of the boundary electrode current from the given current pattern and boundary electrode voltages from the measurements. EIT optimal control problem is fully discretized using the method of finite differences. New Sobolev-Hilbert space is introduced, and the convergence of the sequence of finite-dimensional optimal control problems to EIT coefficient optimal control problem is proved both with respect to functional and control in 2- and 3-dimensional domains.}\\

{{\bf Key words}: Electrical Impedance Tomography, PDE constrained optimal control, method of finite differences, discrete optimal control problem, energy estimate, embedding theorems, convergence in functional, convergence in control}

{{\bf AMS subject classifications}: 35R30, 35Q93, 49M25, 49M05, 65M06, 65N21}




\newpage
{\large 
\section{Introduction}

This paper analyzes inverse EIT problem of estimating an unknown conductivity inside the body based on voltage measurements on the surface of the body when 
electric currents are applied through a set of contact electrodes. Let $Q \in \mathbb{R}^n$ be an open and bounded set  representing body, and assume $ \sigma(x): Q\to \mathbb{R}_+$ be an electrical conductivity function. Electrodes, $(E_l)_{l=1}^m$, with contact impedances vector $Z:= (Z_l)_{l=1}^m \in \mathbb{R}^m_+$ are attached to the periphery of the body, $\partial Q$. Electric current vector $I:=(I_l)_{l=1}^m\in \mathbb{R}^m$ is applied to the electrodes. Vector $I$ is called {\it current pattern} if it satisfies conservation of charge condition
\begin{equation}\label{currentpattern}
\sum_{l=1}^m I_l=0
\end{equation}
The induced constant voltage on electrodes is denoted by $ U:= (U_l)_{l=1}^m \in \mathbb{R}^m$. By specifying ground or zero potential it is assumed that
\begin{equation}\label{grounding}
\sum_{l=1}^m U_l=0
\end{equation}
Let $u: Q\to \mathbb{R}$ is an electrostatic potential. Mathematical model of EIT is described through the following mixed boundary-value problem for second order elliptic PDE:
\begin{alignat}{3}
& -\nabla \cdot (\sigma(x) \nabla u(x)) = 0, & \quad  x\in Q    \label{eq_E:gen_pde} \\
&\sigma(x) \frac{\partial u(x)}{\partial \nu}  = 0,       & \quad  x \in \partial Q -\bigcup\limits_{l=1}^{m} E_{l}\label{eq_E:neumann_zero}\\
& u(x) +  Z_l \sigma(x)\frac{\partial  u(x)}{\partial \nu}  = U_{l} ,       &\quad  x\in E_{l}, \ l= \overline{1,m}\label{eq_E:robin}\\
& \displaystyle\int_{E_l}  \sigma \frac{\partial  u} {\partial \nu} ds = I_l,       &\quad  l= \overline{1,m} \label{eq_E:boundaryflux}
\end{alignat}
where
$\nu$ is the outward normal at $x\in \partial Q$. 
The following is the

{\bf Forward EIT Problem:} {\it Given electrical conductivity map $\sigma$, electrode contact impedance vector $Z$, and electrode {\it current pattern} $I$ it is required to find electrostatic potential $u$ and electrode voltages $U$ satisfying \eqref{eq_E:gen_pde}--\eqref{eq_E:boundaryflux}}.

The goal of the paper is to analyze the following 

{\bf Inverse EIT Problem:} {\it Given electrode contact impedance vector $Z$, electrode {\it current pattern} $I$ and boundary electrode measurement $U^*$, it is required to find electrostatic potential $u$ and electrical conductivity map $\sigma$ satisfying \eqref{eq_E:gen_pde}--\eqref{eq_E:boundaryflux}} with $U=U^*$.

EIT problem has many important applications in medicine, industry, geophysics and material science \cite{holder2004electrical}. We are especially motivated with medical applications on the detection of cancerous tumors from breast tissue or other parts of the body. Relevance of the inverse EIT problem for cancer detection is based on the fact that the conductivity of the cancerous tumor is higher than the conductivity of normal tissues \cite{conductivitytumorhigher}.
Inverse EIT Problem is an ill-posed problem and belongs to the class of so-called Calderon type inverse problems, due to celebrated work \cite{calderon}, where well-posedness of the inverse problem for the identification of the conductivity coefficient $\sigma$ of the second order elliptic PDE \eqref{eq_E:gen_pde}
through Dirichlet-to-Neumann or Neumann-to-Dirichlet boundary maps is presented. 
Significant development in Calderon's inverse problem in the class of smooth conductivity function with spatial dimension $n\geq 3$, concerning questions on uniqueness, stability, reconstruction procedure, reconstruction with partial data was achieved in \cite{kohn1,kohn2,sylvester1987global,nachman1988reconstructions,alessandrini1988stable,astala1,kenig2007,kenig2013}. 
Global uniqueness in spatial dimension $n=2$ and reconstruction procedure through scattering transform and employment of the, so-called $D$-bar method was presented in a key paper \cite{nachman1995}. Further essential development of the $D$-bar method for the reconstruction of discontinuous parameters, regularization due to inaccuracy of measurements, joint recovery of the shape of domain and conductivity are pursued in \cite{dbar1,dbar2,dbar3,dbar4}. 

Mathematical model  \eqref{eq_E:gen_pde}--\eqref{eq_E:boundaryflux} for the EIT, referred to as complete electrode model, was suggested in  \cite{somersalo1992existence}.
This model suggests replacement of the complete potential measurements along the boundary with measurements of constant potential along the electrodes with contact impedances. In \cite{somersalo1992existence} it was demonstrated that the complete electrode model is physically more relevant, and it is capable of predicting the experimentally measured voltages to within 0.1 percent. Existence and uniqueness of the solution to the problem  \eqref{eq_E:gen_pde}-\eqref{eq_E:boundaryflux} was proved in \cite{somersalo1992existence}. {\bf Inverse EIT Problem} is more difficult than the Calderon's problem due to the fact that the infinite-dimensional conductivity function $\sigma$ and finite-dimensional voltage vector $U$ must be identified based on the finitely many boundary electrode voltage measurements.
Hence the input data is finite-dimensional current vector, while in Calderon's problem input data is given via infinite-dimensional boundary operator "Dirichlet-to-Neuman" or "Neuman-to-Dirichlet". \\ Therefore, inverse EIT problem is highly ill-posed and powerful regularization methods are required for its solution. It is essential to note that the size of the input current vector is limited to the number of electrodes, and there is no flexibility to increase its size. It would be natural to suggest that multiple data sets - input currents can be implemented for the identification of the same conductivity function. However, note that besides unknown conductivity function, there is unknown boundary voltage vector with size directly proportional to the size of the input current vector. Accordingly, multiple experiments with "current-to-voltage" measurements is not reducing underdeterminacy of the inverse problem. One can prove uniqueness and stability results by restricting conductivity to the finite-dimensional subset of piecewise analytic functions provided that the number of electrodes is large enough \cite{lechleiter2008,harrach2019}. Within last three decades many methods developed for numerical solution of the ill-posed inverse EIT problem. Without any ambition to present a full review we refer to some significant developments such as recovery of small
inclusions from boundary measurements \cite{ammari1,kwon};  hybrid conductivity imaging
methods \cite{ammari2,seo1,widlak}; multi-frequency EIT imaging methods \cite{ammari3,seo2}; finite element and adaptive finite element method \cite{jin,bangti}; imaging algorithms based on the sparsity reconstruction \cite{ammari3,jin2}; globally convergent method for
shape reconstruction in EIT \cite{harrah};  $D$-bar method, diction reconstruction method, recovering boundary
shape and imaging the anisotropic electrical conductivity \cite{alsaker,dodd,siltanen1,siltanen2,hyvonen}; globally convergent regularization method using Carleman weight function \cite{klibanov}.

Inverse EIT problem was widely studied in the framework of Bayesian statistics \cite{kaipio2005}. In \cite{kaipio2000statistical} inverse EIT problem is formulated as a Bayesian problem of statistical inference and Markov Chain Monte Carlo method with various prior distributions is implemented for calculation of the posterior distributions of the unknown parameters conditioned on measurement data. In \cite{kaipio1999inverse} Bayesian model of the regularized version of the inverse EIT problem is analyzed. In \cite{roininen2014whittle} the Bayesian method with Whittle-Mat\'{e}rn priors is applied to inveres EIT problem. In general the strategy of the Bayesian approach to inverse EIT problem in infinite-dimensional setting is twofold. First approach is based on discretization followed by the application of finite-dimensional Bayesian methods. All the described papers are following this approach, which is outlined in \cite{kaipio2005}. Alternative apprach is based on direct application of the Bayesian methods in functional spaces before discretization \cite{lassas2009,dunlop2015bayesian}.

In this paper we introduce variational formulation of the inverse EIT problem as a PDE constrained coefficient optimal control problem in a new Hilbert space setting. The novelty of the control theoretic model is its adaptation to clinical situation when additional "voltage-to-current" measurements can increase the size of the input data from the number of electrodes $m$ up to $m!$ while keeping the size of the unknown parameters fixed. We pursue discretization of the optimal control problem with the sequence of discrete optimal control problems via finite differences. The main goal of this paper is to prove the convergence of the sequence of finite-dimensional optimal control problems to EIT optimal control problem both with respect to functional and control in 2D and 3D domains. The results on the existence of the optimal control, Fr\'echet differentiability in the Besov space setting, formula for the Fr\'echet gradient, optimality condition, and numerical solution via gradient descent method in 2D model example are addressed in another paper \cite{abdulla2018breast}. 

The organization of the paper is as follows. In Section~\ref{sec:notations} we introduce the notations of the functional spaces. In Section~\ref{sec:Optimal_Control_Problem} we introduce Inverse EIT Problem as PDE constrained optimal control problem.In Section~\ref{sec:discrete_opt_problem} we pursue discretization via finite differences, and introduce approximating sequence of finite dimensional discrete optimal control problems. Section~\ref{sec:Main_Results} formulates the main result. Various key preliminary results are proved in Section~\ref{sec:Preliminary_Results}.
Proof of the main result is completed in Section~\ref{section_E:Approximation Theorem}. Finally, in Section~\ref{sec:conclusions} we outline the main conclusions.

\subsection{Notations}
\label{sec:notations}
Although the main results of the paper are established when number of spatial variables is 2 and 3, for technical reasons we will describe general notations in space of $n$ independent
variables. Differences for the cases $n=2$ or $n=3$ will be specifically mentioned.

Let $Q$ is a bounded domain in $\mathbb{R}^n$; \
$B_r(x)=\{y\in\mathbb R^n: |y-x|<r\}$; \ $m_d(\cdot)$ - $d$-dimensional Lebesgue measure;  We use the standard notation for Banach spaces $C^k(\overline{Q}), \ k\in \mathbb{Z}_*:=\{0\}\cup \mathbb{Z}_+$ of $k$-times continuously differentiable functions on $\overline{Q}$, and we simply write  $C(\overline{Q})$, if $k=0$. The following standard notation will be used for H\"{o}lder spaces:
\begin{itemize}
\item For $k\in \mathbb{Z}_*, 0<\gamma \leq 1$, H\"{o}lder space $C^{k,\gamma}(\overline{Q})$ is the Banach space of elements $u\in C^{k}(\overline{Q})$ with finite norm
\begin{equation*}
\|u\|_{C^{k,\gamma}(\overline{Q})}:=\sum\limits_{|\alpha| \leq k} \|D^\alpha u\|_{C(\overline{Q})}+\sum\limits_{|\alpha| =k}[D^\alpha u]_{C^{0,\gamma}(\overline{Q})}
\end{equation*}
where 
\[ [v]_{C^{0,\gamma}(\overline{Q})}:= \sup\limits_{\substack{x, x'\in Q \\ x\neq x'}} \frac{|v(x)-v(x')|}{|x-x'|^\gamma} \]
\end{itemize}
 Throughout the paper we use standard notations for $L_p(Q), 1\leq p \leq \infty$ spaces; the following standard notations are used for Sobolev spaces \cite{besov79,besov79a}:
\begin{itemize}
\item For $s \in \mathbb{Z}_+, 1\leq p <\infty$,  Sobolev space $W_p^s(Q)$ is the Banach space of measurable functions on $Q$ with finite norm
\begin{equation*}
\| u\|_{W_p^s(Q)} :=
    \left\{
    \begin{array}{l}
    \Big(\int\limits_Q \sum\limits_{|\alpha| \leq s}  |D^{\alpha}u(x)|^p dx \Big)^{\frac{1}{p}}, \quad\text{if} \ 1\leq p <\infty,\\
    \sum\limits_{|\alpha| \leq s} \|D^{\alpha}u(x)\|_{L_\infty(Q)} , \quad\text{if} \  p =\infty,\
    \end{array}\right.
\end{equation*}
where $\alpha = (\alpha_1,..., \alpha_n ) \in \mathbb{Z}_+^n$, $|\alpha| = \alpha_1+...+\alpha_n$, $D_k = \frac{\partial }{\partial x_k}$,  $ D^{\alpha}= D_1^{\alpha_1}...D_n^{\alpha_n}. $
In particular if $p=2$, $H^s(Q) := W_2^s(Q)$  is a Hilbert space with inner product
\begin{equation*}
(f,g)_{H^s(Q)} = \displaystyle\sum_{|\alpha| \leq s}  (D^{\alpha}f(x),D^{\alpha}g(x) )_{L_2(Q)}
\end{equation*}
\item Equivalent inner product and norm in $H^1(Q)$ are given as
\begin{equation*}
((f,g))_{H^1(Q)}:= \int\limits_Q Df\cdot Dg\,dx + \displaystyle\sum\limits_{l=1}^m \int\limits_{E_l} fg\,dS, \ \||f\||_{H^1(Q)}:=((f,f))^{\frac{1}{2}}.
\end{equation*}
\end{itemize}
The following is the new Hilbert space introduced in this paper.
\begin{itemize}
\item $\tilde H^1(Q), \ n=2,3 $ is a linear subspace of $ H^1(Q)$, defined as
\begin{equation*}
    \tilde H^1(Q) = \{ u \in H^1(Q)| u_{x_1 x_2}  \in L_2(Q) \}, \quad \text{if} \ Q\in \mathbb{R}^2
\end{equation*}
\begin{equation*}
    \tilde H^1(Q) = \{ u \in H^1(Q)| u_{x_1 x_2}, u_{x_1 x_3},u_{x_2 x_3},u_{x_1 x_2 x_3}  \in L_2(Q) \}, \quad \text{if} \ Q\in \mathbb{R}^3.
\end{equation*}
$\tilde H_1(Q)$ is an Hilbert space with inner product
\begin{equation*}
    (u,v)_{\tilde H^1(Q)} = 
    \left\{
    \begin{array}{l}
    (u,v)_{H^1}+ (u_{x_1 x_2},v_{x_1 x_2})_{L_2}, \ \text{if} \ n=2\\
    (u,v)_{H^1} + \sum\limits_{\substack{i,j=1\\i<j}}^3  (u_{ x_{i} x_j},v_{ x_{i} x_j})_{L_2} +  (u_{ x_{1} x_2 x_3},v_{ x_{1} x_2 x_3})_{L_2}, \ \text{if} \ n=3
     \end{array}\right.
\end{equation*}
\end{itemize}
Standard notation will be employed for embedding of Banach spaces:
\begin{itemize}
\item $B_1 \hookrightarrow B_2$ means bounded embedding of $B_1$ into $B_2$, i.e. $B_1\subset B_2$, and
\[ \|u\|_{B_2}\leq C \|u\|_{B_1}, \ \forall u\in B_1, \ \text{for some constant} \ C.\]
\item $B_1 \Subset B_2$ denotes compact embedding of $B_1$ into $B_2$, meaning that $B_1 \hookrightarrow B_2$, and every bounded subset of $B_1$ is precompact in $B_2$.
\end{itemize}
\subsection{EIT Optimal Control Problem}
\label{sec:Optimal_Control_Problem}
Consider the optimal control problem on the minimization of the cost functional 
\begin{align}
  \mathcal{J}(v) =   \displaystyle\sum _{l=1}^m \Big |  \displaystyle\int_{E_l} \frac{U_l-u(x)}{Z_l}ds- I_l\Big|^2 + \beta |U-U^*|^2\label{eq_E:cost_functional}
\end{align}
on the control set
\begin{gather} \label{eq_E:control_set}
\mathscr{F}^R=\Big\{v = (\sigma, U)\in \tilde H^1(Q) 
\times \mathbb{R}^m  \,  \Big|  \,  \sum\limits_{l=1}^m U_l = 0, \nonumber
\\
\|\sigma\|^2_{\tilde H^1}+|U|^2 \leq R ^2, \ 
\sigma(x)\geq \sigma_0>0, \, \forall x \in Q \Big \}
\end{gather}
where $\beta> 0, R>0$, and $u=u(\cdot;v)\in H^1(Q)$ is a weak solution of the elliptic problem \eqref{eq_E:gen_pde}--\eqref{eq_E:robin}, i.e.
\begin{equation}\label{eq_E:weak_sol}
\int_Q \sigma \nabla u  \cdot \nabla \eta dx + \displaystyle\sum_{l=1}^m \frac{1}{Z_l}\int_{E_l} u \eta ds = \displaystyle\sum_{l=1}^m \frac{U_l}{Z_l}\int_{E_l} \eta\,ds, \quad  \forall \eta \in H^1(Q).
\end{equation}
This optimal control problem will be called Problem $\mathcal{E}$. 
Note that the first term in the cost functional $\mathcal{J}(v)$ characterizes
the mismatch of the condition \eqref{eq_E:boundaryflux} in light of the Robin condition \eqref{eq_E:robin}.
\subsection{Discrete Optimal Control Problem }
\label{sec:discrete_opt_problem}
To discretize optimal control problems ${\mathscr{E}}$ we pursue finite difference method following the framework introduced in \cite{abdulla2019optimal}. Let $h>0$ and cut $\mathbb{R}^n$ by the planes 
$$x_i = k_i h,   \ i = 1,\ldots,n \ \ \forall\, k_i \in \mathbb{Z}. $$ 
into a collection of elementary cells with length $h$ in each $x_i$-direction. 
For every $h>0$ and  multi-index $\alpha = (k_1,\ldots, k_n)$ we define a cell $C_{h}^{\alpha}$ as
\begin{equation}
C_{h}^{\alpha} = \{ x\in \mathbb{R}^n| \ k_i h \leq x_i \leq (k_i+1) h, \ i = 1,\dots, n      \},
\end{equation}
and consider the collection 
of cells which have non-empty intersection with $Q$:
\begin{equation}
\mathscr{C}_{h}^{Q} = \{ C_{h}^{\alpha} | \ C_{h}^{\alpha} \cap  Q \neq \emptyset  \}
\end{equation}
We now introduce exterior approximation of $\overline{Q}$ as follows:
\begin{equation}
Q_{h} = \bigcup\limits_{C_{h}^{\alpha} \in \mathscr{C}_{h}^Q } C_{h}^{\alpha}
\end{equation}
Obviously, we have $\overline{Q}\subset Q_{h}$. Let $S_{h}=\partial Q_{h}$.
The vertex of the prism $C_{h}^{\alpha}$ whose coordinates are smallest relative to the other vertices, is called its \emph{natural corner}. We are going to identify each prism (cell) by its natural corner. With slight abuse of notation we denote as $C_h^z$, a cell in $\mathbb{R}^n$ of side length $h$ and with natural corner at $z$. Hence, $C_h^\alpha$ and $C_h^{\alpha h}$ are identical.
Consider a lattice 
\[
\mathscr L=\Big\{x\in\mathbb{R}^n~|~\exists\alpha\in\mathbb Z^{n}\text{ s.t. }x_i=k_ih,~i=1,\ldots,n\Big\}.
\]
We will write {$x_{\alpha}=(k_1h,\ldots,k_nh)$}. Bijection $\alpha\mapsto x_{\alpha}$ will henceforth be referred as natural. Given a set $X$ which is in natural bijection with a subset of the set of multi-indexes $\alpha$, we write $\mathscr A(X)$ as the indexing set. Moreover, if $X\subset\mathbb R^n$, then $\mathscr L(X):=\mathscr L\cap X$. When $X=\mathscr L(Y)\subset\mathbb R^n$, we'll agree to write $\mathscr A(Y)$ instead of $\mathscr A(\mathscr L(Y))$. 
These indexes are also in natural bijection with the natural corners of these prisms. In particular, some of the corresponding lattice points may fall on the boundary $S_{h}$. We contrast this set to the set $\mathscr A(Q_{h}')$ of indexes in natural bijection to the lattice points that lie strictly in the interior of $Q_{h}$, and to the set $\mathscr A(Q_{h})$, of all indexes which are in natural bijection with the lattice points that lie in $Q_{h}$. We will write
\[
\sum\limits_{\mathscr A(X)}\quad\text{instead of}\quad\sum\limits_{\alpha\in\mathscr A(X)},
\]
and likewise for other expressions requiring subscripts. We adopt the notation
\[ \alpha\pm e_i := (k_1,...,k_{i}\pm1,...,k_n). \]
To discretize optimal control problem ${\mathscr{E}}$, we need to introduce some refined subsets of grid points of $Q_{h}$.

\begin{equation*}
Q_{h}^{+}=\{x_\alpha\in  Q_{h} :  C_{h}^{\alpha} \cap Q\neq \emptyset \}
\end{equation*}
be a subset of natural corners of the cells in $Q_{h}$. 
We denote as
\begin{equation*}
Q_{h}^{(i)}= \{x_\alpha\in  Q_{h} : x_{\alpha+e_i}\in Q_{h}\}
\end{equation*}
the subset of all grid points $x_\alpha$ in $Q_{h}$ such that the edge $[x_\alpha, x_{\alpha+e_i}] \subset Q_{h}$, and similarly
\begin{equation}
Q_{h}^{(i,j)} = \{ x_{\alpha}\in  Q_{h}: x_{\alpha}+e_i+e_j \in  Q_{h}\}, \quad i,j = 1,2,3 
\end{equation}
Subset of natural corners $x_\alpha$ of cells in $Q_{h}$ which intersect the boundary $S$ is denoted as
\begin{equation}\label{shhat}
\hat{S}_h= \{x_\alpha\in  Q_{h} :  C_{h}^{\alpha} \cap S\neq \emptyset \}
\end{equation}
and
\begin{equation*}
    \hat E_{lh} = \{\ x_\alpha\in   {Q}_{h} : \ E_{l\alpha}:=C_{h}^{\alpha} \cap E_l \neq \emptyset \ \} , \quad l=1, \ldots ,m
\end{equation*}
is a collection of grid pints which are natural corners of $C_{h}^{\alpha} $ containing portion $E_{l\alpha}$ of the boundary curve $E_l$. Let $\Gamma_{l\alpha} = m_{n-1}(E_{l\alpha}), \ l =1, \ldots, m $ is an $n-1$-dimensional Lebesgue measure of $E_{l\alpha}$.
We are going to assume that any control vector $\sigma$ is extended to a larger set $Q+B_1(0)$ as bounded measurable function with preservation of conditions in the control set \eqref{eq_E:control_set}. We introduce discrete grid function by discretizing $\sigma$ through Steklov average: 
\begin{gather}
    \sigma_{\alpha} = \frac{1}{h^n} \int \limits_{x_1}^{x_1+h} \cdots \ \int \limits_{x_n}^{x_n+h} \sigma(y_1, \ldots,y_n)\, dy_1 \cdots \, dy_n, \quad  \alpha\in \mathscr{A}(  Q_{h}), 
    \label{eq_D:steklov_formula}
\end{gather}
where $x_i$ is the i-th coordinate of $ x_{\alpha}$. We use standard notation for finite differences of grid functions $u_\alpha, \sigma_{\alpha} $:
\begin{gather*}
u_{\alpha x_i} = \frac{u_{\alpha+e_i} - u_{\alpha}}{h}, \ u_{\alpha \bar{x_i}} = \frac{u_{\alpha}-u_{\alpha-e_i}}{h},
u_{\alpha x_i} = \frac{u_{\alpha +e_i} -u_{\alpha }}{h}, \quad i = 1,\ldots, n\nonumber\\
\sigma_{\alpha x_i x_j} = \frac{\sigma_{(\alpha+e_j )x_i} -\sigma_{\alpha x_i}}{h} = \frac{\sigma_{\alpha +e_j + e_i}- \sigma_{\alpha+e_j } -  \sigma_{\alpha+e_i } + \sigma_{\alpha }}{h^2},\nonumber\\
\sigma_{\alpha x_1 x_2 x_3} = \frac{\sigma_{(\alpha+e_3 )x_1 x_2} -\sigma_{\alpha x_1 x_2}}{h} = \frac{\sigma_{(\alpha+e_3+e_2 )x_1} -\sigma_{(\alpha+e_3 )x_1} -\sigma_{(\alpha+e_2)x_1 }+ \sigma_{\alpha x_1}}{h^2}=\nonumber\\
\frac{\sigma_{\alpha+e_3+e_2+e_1}-\sigma_{\alpha+e_3+e_2} -\sigma_{\alpha+e_3+e_1}-\sigma_{\alpha+e_2+e_1}+\sigma_{\alpha+e_3}  
+\sigma_{\alpha+e_2}
+\sigma_{\alpha+e_1}-\sigma_{\alpha}}{h^3}.
\end{gather*}
For a given discretization with step size $h$, we employ the notation 
\[ [\eta]_{h} :=\{ \eta_{\alpha}\in \mathbb{R}: \ \alpha \in \mathscr{A}(  Q_{h})\},\]
for the grid function. Next, we define the discrete $\mathscr{H}^1(  Q_{h}), \mathscr{\tilde H}^1(  Q_{h})$ and $L_\infty(Q_h)$ norms:
\begin{gather*}
\| [u]_{h} \|^2_{\mathscr{H}^1(  Q_{h})} :=  \displaystyle\sum_{\mathscr{A}(  Q_{h})}   h^n u^2_{\alpha} + \displaystyle\sum_{i=1}^n \displaystyle \sum_{\mathscr{A}(Q_{h}^{(i)})}   h^n  u^2_{\alpha x_{i}}\nonumber\\
\|| [u]_{h} \||^2_{\mathscr{ H}^1( {Q}_{h})} := \displaystyle\sum_{i=1}^n h^n \displaystyle\sum_{\mathscr{A}(Q_{h}^{(i)})} u_{\alpha x_i}^2 +  \displaystyle\sum_{l=1}^m \displaystyle\sum_{\mathscr{A}(\hat{E}_{l h})} \Gamma_{l \alpha} u_{\alpha}^2\nonumber\\
\| [\sigma]_{h} \|^2_{\mathscr{\tilde H}^1(  Q_{h})} :=  \displaystyle\sum_{\mathscr{A}(  Q_{h})}   h^2 \sigma^2_{\alpha} + \displaystyle\sum_{i=1}^2 \displaystyle \sum_{\mathscr{A}(Q_{h}^{(i)})}   h^2  \sigma^2_{\alpha x_{i}} +  \sum_{\mathscr{A}(  Q_{h}^{+})}  h^2  \sigma^2_{\alpha x_{1} x_2}, \ Q\in \mathbb{R}^2\nonumber\\
\| [\sigma]_{h} \|^2_{\mathscr{\tilde H}^1(  Q_{h})} :=  \displaystyle\sum_{\mathscr{A}(  Q_{h})}   h^3 \sigma^2_{\alpha} + \displaystyle\sum_{i=1}^3 \displaystyle \sum_{\mathscr{A}(Q_{h}^{(i)})}   h^3  \sigma^2_{\alpha x_{i}} + \sum_{\substack{i,j=1\\i<j}}^3\sum_{\mathscr{A} (Q_{h}^{(i,j)})}  h^3  \sigma^2_{\alpha x_{i} x_j}
\nonumber\\
+ \sum_{\mathscr{A} (Q_{h}^{+})}  h^3  \sigma^2_{\alpha x_{1} x_2 x_3}\label{discretenorm3D}, \ Q\in \mathbb{R}^3\nonumber\\
\| [\sigma]_{h} \|_{L_{\infty}(  Q_{h})} :=  \max_{\alpha \in \mathscr{A} (  Q_{h})
} |  \sigma_{ \alpha}|
\end{gather*}
For fixed $R > 0$, define the discrete control sets $\mathscr{F}^R_{h}$ as
\begin{gather}\label{eq_E:disc_control_set}
\mathscr{F}^R_{h} := \bigg\{ [v]_{h} = ([\sigma]_{h}, U) \Big|\, \sum_{l=1}^m U_l = 0, \,  
  \| [\sigma]_{h} \|_{\mathscr{\tilde H}^1(  Q_{h})}^2 + |U|_{\mathbb{R}^m}^2 \leq R^2,
\nonumber\\
\sigma_{ \alpha}\geq \sigma_0>0, \, \forall \alpha \in \mathscr{A} (  Q_{h})
\bigg\}
\end{gather} 
and the interpolating map $\mathscr{P}_{h}$ as
\begin{gather}
\mathscr{P}_{h} : \bigcup \limits_{R}\mathscr{F}^R_{h} \to \bigcup \limits_{R}\mathscr{F}^R, \qquad \mathscr{P}_{h}([v]_{h}) = (\mathscr{P}_{h}([\sigma]_{h}), U) =  (\sigma_{h}',U)
\nonumber
\end{gather}
where $\sigma_{h}'$ is a multilinear interpolation of $[\sigma]_h$, which assigns the value $\sigma_{\alpha}$ to each grid point of $C_{h}^{\alpha}$, and it is a piecewise linear with respect to each variable $x_i$ when the other variables are fixed. Precisely,
\begin{gather}
 \sigma_{h}'(x) = \sigma_{\alpha} + \sigma_{\alpha x_1}(x_1-k_1 h)  + \sigma_{\alpha x_2}(x_2-k_2 h)\nonumber\\ +  \sigma_{\alpha x_1 x_2}(x_1-k_1 h) (x_2-k_2 h), x\in C_h^\alpha, \ \ n=2\\
    \sigma_{h}'(x) = \sigma_{\alpha} + \displaystyle\sum_{i=1}^3 \sigma_{\alpha x_i}(x_i-k_i h)
    + \sum_{\substack{i,j=1\\i<j}}^3 \sigma_{\alpha x_i x_j}(x_i-k_i h) (x_j-k_j h) \nonumber
    \\+ \sigma_{\alpha x_1 x_2 x_3} \prod \limits_{1\leq i \leq 3} (x_i-k_i h), x\in C_h^\alpha, \ \ n=3.  \label{eq_E_3d:multilinear_formula}
\end{gather}
We also define the discretizing map $\mathscr{Q}_{h}$ as
\begin{gather}
\mathscr{Q}_{h} : \bigcup \limits_{R}\mathscr{F}^R \to \bigcup \limits_{R}\mathscr{F}^R_{h}, \qquad \mathscr{Q}_{h}( v) = ( \mathscr{Q}_{h}( \sigma),U) =([\sigma]_{h},U)
\nonumber
\end{gather}
where $[\sigma]_{h}  =\{ \sigma_{\alpha}\} $, with $\sigma_{\alpha}$ given by \eqref{eq_D:steklov_formula} for each $\alpha \in \mathscr{A}(  Q_{h})$. 

Next, we define a discrete state vector, which is a solution of the discretized elliptic problem \eqref{eq_E:gen_pde}--\eqref{eq_E:robin}.
\begin{definition} \label{dfn_E:disc_weak_sln}
Given $[v]_{h}$, the grid function $[u([v]_{h})  ]_{h}$ is called a discrete state vector of problem $\mathcal E$ if it satisfies 
\begin{gather}
h^n \displaystyle\sum_{\mathscr{A}(Q_{h}^{+})} \sigma_{\alpha} \displaystyle\sum_{i=1}^n  u_{\alpha x_i} \eta_{\alpha x_i}  + \displaystyle\sum_{l=1}^m \frac{1}{Z_l}\displaystyle\sum_{\mathscr{A}(\hat{E}_{l h})} \Gamma_{l \alpha} u_{\alpha}   \eta_{\alpha} + J_h([u]_h, [\eta]_h)\nonumber\\
 = \displaystyle\sum_{l=1}^m \frac{U_l}{Z_l}\displaystyle\sum_{\mathscr{A}(\hat{E}_{lh})} \Gamma_{l \alpha} \eta_{\alpha}
\label{eq_E:disc_weak_sln}
\end{gather}
\end{definition}
for arbitrary grid function $[\eta]_h$, where
\begin{equation*}
    J_h([u]_h, [\eta]_h) = h^n \displaystyle\sum_{\mathscr{A}({S}_{h})}  \displaystyle\sum_{i=1}^n \theta_{\alpha }^i u_{{\alpha} x_i} \eta_{{\alpha} x_i}, \     \theta_{\alpha}^i = \left\{
    \begin{array}{cl}
    1  &{\rm if} \ \alpha \in \mathscr{A}(Q_{h}^{(i)} \setminus Q_{h}^{+})    \\
    0 & {\rm otherwise}
    \end{array}
    \right.
\end{equation*}
The necessity of adding $J_{\alpha}$ to \eqref{eq_E:disc_weak_sln} is that some $u_{\alpha x_i}$ and $\eta_{\alpha x_i}$ values on $S_{h}$ are not present in the first term 
of \eqref{eq_E:disc_weak_sln}. Addition of these terms to \eqref{eq_E:disc_weak_sln} through 
$J_{\alpha}$ is essential for the proof of stability of our discrete scheme.

In Section \ref{sec:Preliminary_Results}, it will be proved that for a given $[v]_{h}\in \mathscr{F}^R_{h}$ there exists a unique discrete state vector of problem $\mathcal E$. Consider minimization of the discrete cost functional
\begin{equation}\label{eq_E:discerte_obj_func}
    \mathcal{J}_{h}( [v]_{h} ) =  \displaystyle\sum_{l=1}^m  \Big( \displaystyle\sum_{ \mathscr{A}(\hat E_{lh})} \Gamma_{l \alpha}\frac{U_l-u_{\alpha}}{Z_l}   - I_l \Big)^2 + \beta |U-U^*|^2 
\end{equation} 
on a control set $\mathscr{F}^R_{h}$, where $u_\alpha$'s are components of the discrete state vector $[u([v]_{h})]_{h}$ of the Problem $\mathscr{E}$. The formulated discrete optimal control problem will be called Problem $\mathscr{E}_{h}$.

Next, we define three interpolations of the discrete state vector $[u]_h$. Piecewise constant interpolation $\tilde u_{h}: Q_h\to \mathbb{R}$ assigns
to the interior of each cell in $Q_{h}$ the value of $u_{\alpha}$ at its natural corner:
\begin{equation}
 \tilde u_{h} \big|_{C_{h}^{\alpha '}} =u_{\alpha}, \qquad \forall \alpha \in \mathscr{A} (Q_{h}^+). 
\end{equation}
Piecewise constant interpolation of the discrete $x_i$-derivative $\tilde u_{h}^i : Q_h \to \mathbb{R}, \ i= 1,..,n$ assign to the interior of each cell in $Q_{h}$ the value of the forward spatial difference at the natural corner:
\begin{equation}
 \tilde u_{h}^i \big|_{C_{h}^{\alpha '}} =u_{\alpha x_i}, \qquad \forall \alpha \in \mathscr{A} (Q_{h}^+). 
\end{equation}
Multilinear interpolation  $u'_{h} : Q_h \to \mathbb{R}$ assigns the value $u_{\alpha}$ to each grid point in $\mathscr{L} ({Q}_{h})$, and it is a peicewise linear with respect to each variable $x_i$ when the rest of variables are fixed. 
\section{Main Result}
The following is the main result of the paper.
\label{sec:Main_Results}
\begin{theorem} \label{thm_E:convergence-of-optimal-control}
The sequence of discrete optimal control problems $\mathcal{E}_{h}$ approximates the optimal control problem $\mathcal{E}$ with respect to functional, i.e.
\begin{equation}
    \lim_{h \to 0} \mathscr{J}_{h_*} = \mathscr{J}_{*},\label{eq_E:conv_of_J}
\end{equation}
where 
\begin{equation}
    \mathscr{J}_{h_*} = \inf_{ \mathscr{F}^R_{h}} \mathscr{J}_{h} ([v]_{h}), \ \mathscr{J}_*:=\inf\limits_{f\in\mathscr{F}^R}\mathscr{J}(v).
\end{equation}
Furthermore, let $\{\epsilon_{h}\}$ be a sequence of positive real numbers with $\displaystyle\lim_ {h  \to 0}  \epsilon_{h} = 0$. If the sequence $[v]_{h,\epsilon}=([\sigma]_{h,\epsilon},U^{h,\epsilon}) \in  \mathscr{F}^R_{h} $ is chosen so that
\begin{equation}
    \mathscr{J}_{h_*} \leq  \mathscr{J}_{h} ([v]_{h,\epsilon }) \leq \mathscr{J}_{h_*} + \epsilon_{h},
\end{equation}
then we have
\begin{equation}
    \displaystyle\lim_ {h  \to 0} \mathscr{J}(\mathscr{P}_{h} ([v]_{h,\epsilon }) ) = \mathscr{J}_{*},\label{eq_E:conv_of_J_over_interp}
\end{equation}
the sequence $\{ (\mathscr{P}_{h} ([\sigma]_{h,\epsilon }), U^{h,\epsilon}) \}$ 
\begin{itemize}
\item is precompact in Tikhonov topology of $\tilde H^1(Q)\times \mathbb{R}^m$ formed with the product of the weak topology of $\tilde H^1(Q)$ and Euclidean topology of $\mathbb{R}^m$;
\item is precompact in Tikhonov topology of $C^{0,\mu}(\overline{Q})\times \mathbb{R}^m$, $0<\mu < \frac{1}{2}$ formed with the product of the strong topology of H\"{o}lder space $C^{0,\mu}(\overline{Q})$ and Euclidean topology of $\mathbb{R}^m$;
\end{itemize}
and all the corresponding limit points $v_*=(\sigma_*,U_*)$ are optimal controls of the problem $\mathcal{E}$. Moreover, if $v_*=(\sigma_*,U_*)$ is any such limit point, then there exists a subsequence $h'$ such that the multilinear interpolations $u'_{h'}$ of the discrete state vectors $[u([v]_{h',\epsilon })]_{h'}$  converge to weak solution $u = u(x; v_*)$ of the elliptic problem  \eqref{eq_E:gen_pde}-\eqref{eq_E:robin}, weakly in $H^1(Q)$, strongly in $L_2(Q)$, and almost everywhere on $Q$. 
\end{theorem}
\section{Preliminary Results}
\label{sec:Preliminary_Results}
The following lemma presents a key discrete energy estimate for the elliptic PDE problem:
\begin{lemma}[Discrete Energy Estimate] \label{lem_E:Discrete_Energy_Estimate} 
For any $[v]_h\in \mathscr{F}^R_{h}$, discrete state vector $[u([v]_{h})  ]_{h}$ satisfies the energy estimate:
\begin{equation}
\|| [u]_{h} \||_{\mathscr{ H}^1( {Q}_{h})} \leq M |U|,
\label{eq_E:disc_energy_est}
\end{equation}
where $M$ depends on $\sigma_0, Z$ and $Q$.
\end{lemma}
\textbf{Proof:} 
We set $ \eta_{\alpha} = u_{\alpha}$ in \eqref{eq_E:disc_weak_sln} to get
\begin{equation*}
h^n \displaystyle\sum_{\mathscr{A}(Q_{h}^{+})} \sigma_{\alpha} \displaystyle\sum_{i=1}^n  u_{\alpha x_i}^2   + \displaystyle\sum_{l=1}^m \frac{1}{Z_l}\displaystyle\sum_{\mathscr{A}(\hat{E}_{l h})} \Gamma_{l \alpha} u_{\alpha}^2 + J_{\alpha}(u_{\alpha}, u_{\alpha}) = \displaystyle\sum_{l=1}^m \frac{U_l}{Z_l}\displaystyle\sum_{\mathscr{A}(\hat{E}_{l h})} \Gamma_{l \alpha}u_{\alpha},
\end{equation*}
and by recalling the definition of $J_{\alpha}$ and the fact that $0 <\sigma_0\leq \sigma_{ \alpha}$ we have 
\begin{equation} \label{eq_E:disc_energy_in_progress1}
\mu \|| [u([v]_{h})  ]_{h}\||_{\mathscr{ H}^1( {Q}_{h})}^2 \leq  \displaystyle\sum_{l=1}^m Z_l^{-1}U_l\displaystyle\sum_{\mathscr{A}(\hat{E}_{lh})}\Gamma_{l \alpha} u_{\alpha}
\end{equation}
where $\mu = \min\{\sigma_0 , \min \limits_{l}Z_l^{-1} \} $. Using Cauchy-Schwarz inequality we derive 
\begin{gather}
      \displaystyle\sum_{l=1}^m Z_l^{-1}U_l\displaystyle\sum_{\mathscr{A}(\hat{E}_{lh})}\Gamma_{l \alpha} u_{\alpha} \leq 
        m_{n-1}^{\frac{1}{2}}(\partial Q) \max\limits_{l}Z_l^{-1} |U| \  \|| [u([v]_{h})  ]_{h}\||_{\mathscr{ H}^1( {Q}_{h})}. \label{est1}
\end{gather}
From \eqref{eq_E:disc_energy_in_progress1} and \eqref{est1}, \eqref{eq_E:disc_energy_est} follows with $M=\mu^{-1} m_{n-1}^{\frac{1}{2}}(\partial Q) \max\limits_{l}Z_l^{-1}$. $\quad\square$
\begin{corollary} \label{cor_E:sys_for_disc_state} For any $[v]_{h}\in \mathscr{F}_{h}^R$, there exists a unique discrete state vector $[u([v]_{h})]_{h}$.
\end{corollary}
Assertion of the corollary follows from energy estimate with similar arguments as in \cite{ladyzhenskaya2013boundary}. By replacing $u_{x_i}$, $\eta_{x_i}$ with respective difference quotients, from \eqref{eq_E:disc_weak_sln} it follows
\begin{gather}
\sum_{\mathscr{A}(Q_h)} \Big \{ \mathscr{L}_\alpha \cdot [u]_h -\mathscr{G}_\alpha(U)\Big \} \eta_\alpha =0,\label{SLAE1}
\end{gather}
where $\mathscr{L}_\alpha$ is a vector of the same size as $[u]_h$ and $\mathscr{G}_\alpha: \mathbb{R}^m\to \mathbb{R}$ is a linear functional. Since the values of $\eta_\alpha$ are independent, \eqref{SLAE1} is equivalent to the following system of linear algebraic equations (SLAE)
\begin{gather}
\mathscr{L}_\alpha \cdot [u]_h =\mathscr{G}_\alpha(U), \ \alpha \in {\mathscr{A}(Q_h)}. \label{SLAE2}
\end{gather}
Note that the number of equations, and the number of unknowns $u_\alpha$ in \eqref{SLAE2} are both equal to number of vertices in a grid $Q_h$. Addition of the expression $J_\alpha(u_\alpha,\eta_\alpha)$ to the discrete identity \eqref{eq_E:disc_weak_sln} surved exactly to this purpose. From the energy estimate \eqref{eq_E:disc_energy_est} it easily follows that the corresponding homogeneous SLAE has only a zero solution. Therefore, claim of the corollary is a consequence of the well-known result of linear algebra.

Another crucial consequence of the energy estimate \eqref{eq_E:disc_energy_est} is uniform $H^1(Q)$-bounded of the interpolations of the discrete state vector:  
\begin{corollary}\label{H1boundforu} Multilinear interpolation $u_h'$ of the discrete state vector is uniformly bounded in $H^1(Q)$:
\begin{equation}\label{H1energyestimate}
\sup\limits_{[v]_h\in \mathscr{F}_{h}^R} \|u_h'\|_{H^1(Q)}\leq C,
\end{equation}
where $C$ depends on $\sigma_0, Z, Q, R, n$. 
\end{corollary} 
Indeed, first of all from \cite{abdulla2019optimal} (formula (4.13)) it follows that 
\begin{equation}
 \int_{Q} |Du'_{h}|^2 dx 
\leq 2^{n-1} \displaystyle\sum_{i=1}^n\displaystyle\sum_{\mathscr{A}(Q_h^+)} h^n  |u_{\alpha x_i }|^2.
\label{eq_E:l2_norm_upper_interp}
\end{equation}
Next, we establish that the sequences $u_h^{\prime}$ and $\tilde u_{h}$ are equivalent in strong topology of $L_2(E_l)$, as $h\to 0$. The proof is similar to the statement (d) of Theorem 14 in \cite{abdulla2019optimal}. The following estimate is proved in \cite{abdulla2019optimal} (estimate (4.23)):
\begin{gather}
 |\tilde u_{h}(x)-u_h^{\prime}(x)|\leq (2^n-1)n \sum\limits_{{\text{edges of} \ C_{h}^{\alpha}}} h^2 |u_{\alpha'x_j}|^2, \ x\in C_{h}^{\alpha}, \label{L2Sequiv}
 \end{gather}
where the summation on the right-hand side is taken over all $\alpha'$ and $j$ such that $\alpha'\in \mathscr{A}(C_{h}^{\alpha})$ and $\alpha'+e_j\in \mathscr{A}(C_{h}^{\alpha})$.
Since, $E_l$ is Lipschitz, we have $m_{n-1}(E_{l\alpha})\leq L h^{n-1}$, where $L$ is a Lipschitz constant of the boundary $S$. Therefore, from \eqref{L2Sequiv} and \eqref{eq_E:disc_energy_est} it follows that
\begin{gather}\label{L2equiv1}
\|\tilde u_{h}-u_h^{\prime}\|^2_{L_2(E_l)}\leq L(2^n-1) 2^{n-1}n \displaystyle\sum_{i=1}^n \displaystyle \sum_{\mathscr{A}(Q_{h}^{(i)})}   h^{n+1}  u^2_{\alpha x_{i}} \to 0, 
\end{gather}
as $h\to 0$. Assuming that $h\leq 1$, from \eqref{eq_E:l2_norm_upper_interp} and \eqref{L2equiv1} it follows that
\begin{equation}\label{H1boundtriple}
\||u_h'\||_{H^1(Q)}^2 \leq C_1 \||[u]_h\||^2_{\mathscr{ H}^1( {Q}_{h})}, 
\end{equation}
where $C_1=2^{n-1}(2L(2^n-1)nm+1)$. Due to equivalency of the norms $\|\cdot\|$ and $\||\cdot\||$ in $H^1(Q)$ (see Lemma 5.1 in \cite{abdulla2018breast}), from \eqref{H1boundtriple}, \eqref{H1energyestimate} follows.

Discrete energy estimate implies the following interploation
\begin{lemma}\label{thm_E:interp_relations}
Let $R>0$ is fixed, and $\{ [v]_{h} \}$ is a sequence of discrete control vectors such that $[v]_h \in \mathscr{F}_{h}^R$ for each $h$. Then the following statements hold:

(\textit{a}) The sequences $\{u'_{h} \}$ and $\{\tilde u_{h} \}$ are uniformly bounded in $L_{2}(Q_{h})$.

(\textit{b}) For each $i \in \{1,\ldots,n\}$, the sequences $ \{\tilde u_{h}^i \},\, \{ \frac{\partial u'_{h}}{\partial x_i}  \}$ are uniformly bounded
in $L_{2}(Q_{h})$. 

(\textit{c}) the sequence $\{ \tilde u_{h} - u'_{h} \}$ converges strongly to 0 in $L_{2}(Q)$ as $h \to 0$.

(\textit{d}) For each $i \in \{1,\ldots,n\}$, the sequences $\{ \frac{\partial u'_{h}}{\partial x_i} - \tilde u_{h}^i \}$ converges weakly to zero in $L_{2}(Q)$ as $h \to 0$.

(\textit{e}) the sequence $\{ \tilde u_{h} - u'_{h} \}$ converges strongly to 0 in $L_{2}(S)$ as $h \to 0$.
\end{lemma}
The proof of the claims (a)-(d) coincides with the proofs of similar claims in Theorem 14 of \cite{abdulla2019optimal}. The claim (e) is proved above in \eqref{L2equiv1}.

Next, we recall the necessary and sufficient condition for the convergence of the discrete optimal control problems $\mathscr{E}_h$, which is the suitable criteria to employ for the proof of method of finite differences for the optimal control problems with distributed parameters (\cite{abdulla13}-\cite{abdulla-cosgrove2020}). 
\begin{lemma} \label{lem_E:necessary condition} \cite{vasil1981methods} The sequence of discrete optimal control problems $\mathscr{E}_{h}$ approximates the continuous optimal control problem $\mathscr{E}$ with respect to the functional if and only if the following conditions are satisfied:
\begin{enumerate}
\item  For arbitrary sufficiently small $\epsilon>0$ there exists $h_1 =h_1(\epsilon)$ such that $\mathscr{Q}_{h}(v) \in \mathscr F^R_{h}$
for all $v \in \mathscr F^{R-\epsilon}$ and $h\leq h_1$; Moreover, for any fixed $\epsilon>0$ and for all $v \in \mathscr F^{R-\epsilon}$ the following inequality is satisfied:
\begin{equation}\label{criteriacond1}
    \lim\sup_{h \to 0}(\mathscr{J}_{h}( \mathscr{Q}_{h}(v)) - \mathscr{J}(v) ) \leq 0.
\end{equation}
\item  For arbitrary sufficiently small $\epsilon>0$ there exists $h_2 =h_2(\epsilon)$ such that $\mathscr{P}_{h}([v]_{h}) \in \mathscr F^{R+\epsilon}$ for all $[v]_{h} \in \mathscr F^R_{h}$ and $h\leq h_2$;  moreover, for all $[v]_{h} \in \mathscr F^R_{h}$,  the following inequality is satisfied:
\begin{equation}\label{criteriacond2}
    \lim\sup_{h \to 0}(\mathscr{J}( \mathscr{P}_{h}([v]_{h})) - \mathscr{J}_{h}([v]_{h}) ) \leq 0.
\end{equation}
\item  For arbitrary sufficiently small $\epsilon>0$, the following inequalities are satisfied:
\begin{equation*}
\lim\sup_{\epsilon \to 0} \mathcal{J}_*(\epsilon ) \geq \mathcal{J}_*, \  \lim\inf_{\epsilon \to 0} \mathcal{J}_*(-\epsilon ) \leq \mathcal{J}_*,
\end{equation*}
where $\mathcal{J}_*(\pm \epsilon) = \displaystyle\inf_{\mathscr F^{R\pm \epsilon}} \mathcal{J}(v)$.
\end{enumerate}
\end{lemma}
Our next goal is to show that the mappings $\mathscr{P}_{h}$ and $\mathscr{Q}_{h}$ satisfy the conditions of Lemma \ref{lem_E:necessary condition}. The following lemma plays a key role to prove this claim. The proof is similar to the proof of Proposition 11 in \cite{abdulla2019optimal}.  
\begin{lemma}  \label{lemma_E:finite-diff_bound} Let $Q \in \mathbb{R}^2$. Then for $\forall \epsilon >0$, there exists $\delta>0$ such that 
$$\displaystyle\sum_{\mathscr{A}(Q_{h}^{+})} h^2 |\sigma_{\alpha x_1 x_2}|^2 \leq (1+\epsilon) \Big\|\frac{\partial^2 \sigma}{\partial x_1 \partial x_2} \Big\|^2_{L_2(Q_{h})}$$
whenever $h<\delta$.
\end{lemma}
\textit{Proof:}
For each $h>0$, define the function $\tilde \sigma_h^{12}$ as 
\begin{equation}
    \tilde \sigma_h^{12} \Big|_{C_{h}^{\alpha}} =  \sigma_{\alpha x_1 x_2}, \quad \forall \alpha \in \mathscr{A}(Q_{h}^{+})
\end{equation}
In the following we will prove that 
\begin{equation}
    \tilde \sigma_h^{12}  \to \frac{\partial^2 \sigma}{\partial x_1 \partial x_2}  \quad \text{strongly in $L_2(Q)$  as } h \to 0
\end{equation}
As an element of $\tilde H^1(Q)$, almost all restrictions of $\sigma$ to lines parallel to the $x_i$ direction are
absolutely continuous, moreover, restrictions of $\frac{\partial \sigma}{\partial x_1}$ (or $\frac{\partial \sigma}{\partial x_2}$) to lines parallel to the $x_2$ (or $x_1$) direction are
absolutely continuous.
Therefore, for almost every $z=(z_1,z_2)\in Q$ we have
\begin{gather}
   \int\limits_{C_h^z}\frac{\partial^2 \sigma}{\partial x_1 \partial x_2}\,dy
    = \sigma(z+he_2+he_1) -\sigma(z+he_2) -\sigma(z+he_1)+\sigma(z)\label{2nddiffint}
\end{gather}
In the following transformation, we write simply $\mathscr{A}$ instead of summation index set $\mathscr{A}(Q^{+}_{h})$. Using the definition of Steklov average \eqref{eq_D:steklov_formula}, \eqref{2nddiffint} and Cauchy-Schwartz inequality, we get
\begin{gather}
   \Big \| \tilde \sigma_h^{12} - \frac{\partial^2 \sigma}{\partial x_1 \partial x_2} \Big\|^2_{L_2(Q_{h})} = 
\displaystyle\sum_{ \mathscr{A}} \int\limits_{C_{h}^{\alpha}} \Big|\sigma_{\alpha x_1 x_2} - \frac{\partial^2 \sigma(x)}{\partial x_1 \partial x_2}\Big|^2 dx=
    \nonumber
   \\ \displaystyle\sum_{ \mathscr{A}}  \int\limits_{C_{h}^{\alpha}} \Big| \frac{1}{h^4}\Big[ \int\limits_{C_{h}^{\alpha+e_1+e_2}}
   dz -\int\limits_{C_{h}^{\alpha+e_1}}dz 
   -\int\limits_{C_{h}^{\alpha+e_2}}dz 
      +\int\limits_{C_{h}^{\alpha}}dz \sigma(z) \Big]- \frac{\partial^2 \sigma(x)}{\partial x_1 \partial x_2} \Big|^2 dx=
    \nonumber
   \\ \displaystyle\sum_{ \mathscr{A}}  \int\limits_{C_{h}^{\alpha}} \Big| \frac{1}{h^4} \int\limits_{C_{h}^{\alpha}} \Big [\sigma |_{z+he_1+he_2} 
   -\sigma |_{z+he_2} -\sigma |_{z+he_1}+\sigma |_z\Big]\,dz
   - \frac{\partial^2 \sigma(x)}{\partial x_1 \partial x_2} \Big|^2 dx\nonumber\\
   = \displaystyle\sum_{ \mathscr{A}} \frac{1}{h^8} \int\limits_{C_{h}^{\alpha}} \Big|  \int\limits_{C_{h}^{\alpha}}  \int\limits_{C_h^z}\Big [\frac{\partial^2 \sigma(y)}{\partial x_1 \partial x_2} - \frac{\partial^2 \sigma(x)}{\partial x_1 \partial x_2}\Big ] \,dy\,dz
    \Big|^2 dx
         \nonumber
        \\
    \leq  \displaystyle\sum_{ \mathscr{A}} \frac{1}{h^4} \int\limits_{C_{h}^{\alpha}}   \int\limits_{C_{h}^{\alpha}}  \int\limits_{C_h^z} \Big| \frac{\partial^2 \sigma(y)}{\partial x_1 \partial x_2} - \frac{\partial^2 \sigma(x)}{\partial x_1 \partial x_2}\Big|^2\,dy\,dz
     dx\label{steklovestimate}
\end{gather}
Changing integration order with respect to $y$ and $z$, we have
\begin{gather}
 \int\limits_{C_{h}^{\alpha}}  \int\limits_{C_h^z} \Big| \frac{\partial^2 \sigma(y)}{\partial x_1 \partial x_2} - \frac{\partial^2 \sigma(x)}{\partial x_1 \partial x_2}\Big|^2\,dy\,dz= \Big (\int\limits_{C_{h}^{\alpha}}(y_1-k_1h)(y_2-k_2h)\,dy+\nonumber\\
\int\limits_{C_{h}^{\alpha+e_1}}((k_1+2)h-y_1)(y_2-k_2h)\,dy + \int\limits_{C_{h}^{\alpha+e_2}}(y_1-k_1h)((k_2+2)h-y_2)\,dy\nonumber\\
\int\limits_{C_{h}^{\alpha+e_1+e_2}}((k_1+2)h-y_1)((k_2+2)h-y_2)\,dy\Big )  \Big| \frac{\partial^2 \sigma(y)}{\partial x_1 \partial x_2} - \frac{\partial^2 \sigma(x)}{\partial x_1 \partial x_2}\Big|^2\leq h^2 \times \nonumber\\
\Big ( \int\limits_{C_{h}^{\alpha}}\,dy+ \int\limits_{C_{h}^{\alpha+e_1}}\,dy+\int\limits_{C_{h}^{\alpha+e_2}}\,dy+\int\limits_{C_{h}^{\alpha+e_1+e_2}}\,dy\Big ) \Big| \frac{\partial^2 \sigma(y)}{\partial x_1 \partial x_2} - \frac{\partial^2 \sigma(x)}{\partial x_1 \partial x_2}\Big|^2 \label{fubini1}
\end{gather}
From \eqref{steklovestimate}, \eqref{fubini1} it follows that
\begin{gather}
   \Big \| \tilde \sigma_h^{12} - \frac{\partial^2 \sigma}{\partial x_1 \partial x_2} \Big\|^2_{L_2(Q_{h})} 
   \leq \displaystyle\sum_{ \mathscr{A}} \frac{1}{h^2}\int\limits_{C_{h}^{\alpha}} \Big(\int\limits_{C_{h}^{\alpha+e_1+e_2}}
   dz
   \nonumber
    \\
    +\int\limits_{C_{h}^{\alpha+e_1}}dz+\int\limits_{C_{h}^{\alpha+e_2}}dz+\int\limits_{C_{h}^{\alpha}}dz\Big ) \Big| \frac{\partial^2 \sigma(z)}{\partial x_1 \partial x_2}  - \frac{\partial^2 \sigma(x)}{\partial x_1 \partial x_2} \Big|^2\,dx \label{eq_E:four_terms}
\end{gather}
Let  $\forall \ \epsilon > 0$ is fixed. Since $C^2 (\overline{Q+B_1(0)})$ is dense in $\tilde H^1(Q+B_1(0))$ we can choose $g\in C^2 (\overline{Q+B_1(0)})$ such that
\begin{equation}\label{L2approx2der}
    \Big \|\frac{\partial^2 \sigma}{\partial x_1 \partial x_2}-\frac{\partial^2 g}{\partial x_1 \partial x_2}\Big \|_{L_2(Q+B_1(0))}^2 < \frac{\epsilon}{24(1+m_n(Q))}.
\end{equation}
From \eqref{eq_E:four_terms} it follows 
\begin{gather}
   \Big \| \tilde \sigma_h^{12} - \frac{\partial^2 \sigma}{\partial x_1 \partial x_2} \Big\|^2_{L_2(Q_{h})} 
   \leq \displaystyle\sum_{ \mathscr{A}} \frac{3}{h^2}\int\limits_{C_{h}^{\alpha}} (I_1 + I_2 +I_3)\,dx,\label{epsilonbound}
\end{gather}
where 
\begin{gather*}
    I_1 = \Big(\int\limits_{C_{h}^{\alpha+e_1+e_2}}
   dz +
   \int\limits_{C_{h}^{\alpha+e_1}}dz+\int\limits_{C_{h}^{\alpha+e_2}}dz+\int\limits_{C_{h}^{\alpha}}\,dz\Big ) \Big| \frac{\partial^2 \sigma (z) }{\partial x_1 \partial x_2}  -\frac{\partial^2 g (z)}{\partial x_1 \partial x_2} \Big|^2,\nonumber\\  
I_2 =  \Big(\int\limits_{C_{h}^{\alpha+e_1+e_2}}
   dz +
   \int\limits_{C_{h}^{\alpha+e_1}}dz+\int\limits_{C_{h}^{\alpha+e_2}}dz+\int\limits_{C_{h}^{\alpha}}\,dz\Big ) \Big| \frac{\partial^2 g (z)}{\partial x_1 \partial x_2}  - \frac{\partial^2 g (x)}{\partial x_1 \partial x_2} \Big|^2,\nonumber\\
    I_3 =  \Big(\int\limits_{C_{h}^{\alpha+e_1+e_2}}
   dz +
   \int\limits_{C_{h}^{\alpha+e_1}}dz+\int\limits_{C_{h}^{\alpha+e_2}}dz+\int\limits_{C_{h}^{\alpha}}\,dz\Big ) \Big| \frac{\partial^2 g (x)}{\partial x_1 \partial x_2} - \frac{\partial^2 \sigma (x)}{\partial x_1 \partial x_2} \Big|^2.
\end{gather*}
Since $\frac{\partial^2 g }{\partial x_1 \partial x_2}$ is uniformly continuous on $Q+B_1(0)$, there exists $\delta = \delta(\epsilon) > 0$ such that
\begin{equation}\label{unifcont}
    \Big |\frac{\partial^2 g (z) }{\partial x_1 \partial x_2}- \frac{\partial^2 g (x) }{\partial x_1 \partial x_2} \Big |^2 < \frac{\epsilon}{24(1+m_n(Q))}
\end{equation}
whenever $|z - x| < \delta $. Let $h_{\epsilon} > 0$ satisfy
\begin{equation}
    \sqrt{8}\, h_{\epsilon} < \delta, 
\end{equation}
Then \eqref{unifcont} is satisfied for each  $h < h_{\epsilon}$, any $\alpha \in \mathscr{A}$ , and any $x, z \in C_{h}^{\alpha+e_1+e_2} \cup C_{h}^{\alpha+e_1} \cup C_{h}^{\alpha+e_2}\cup C_{h}^{\alpha}$. Assuming $h$ is chosen so small that $m_n(Q_h)\leq 2 m_n(Q)$, from \eqref{L2approx2der}, \eqref{unifcont} it follows
\begin{gather*}
    \displaystyle\sum_{ \mathscr{A}} \frac{3}{h^2}\int\limits_{C_{h}^{\alpha}} I_1\,dx \leq 12  \Big \|\frac{\partial^2 \sigma}{\partial x_1 \partial x_2}-\frac{\partial^2 g}{\partial x_1 \partial x_2}\Big \|_{L_2(Q+B_1(0))}^2  <  \frac{\epsilon}{2(1+m_n(Q))}\\
 \displaystyle\sum_{ \mathscr{A}} \frac{3}{h^2}\int\limits_{C_{h}^{\alpha}} I_2\,dx <   \frac{\epsilon m_n(Q)}{1+m_n(Q)}\\
   \displaystyle\sum_{ \mathscr{A}} \frac{3}{h^2}\int\limits_{C_{h}^{\alpha}} I_3\,dx \leq 12  \Big \|\frac{\partial^2 \sigma}{\partial x_1 \partial x_2}-\frac{\partial^2 g}{\partial x_1 \partial x_2}\Big \|_{L_2(Q+B_1(0))}^2  <  \frac{\epsilon}{2(1+m_n(Q))}
\end{gather*}
From \eqref{epsilonbound} we deduce
\begin{equation}
   \Big \| \tilde \sigma_h^{12} - \frac{\partial^2 \sigma}{\partial x_1 \partial x_2} \Big\|^2_{L_2(Q_{h})} 
   <  \epsilon, \quad \forall h \leq h_{\epsilon}
\end{equation}
Lemma is proved. $\square$

The following lemma expresses similar result for 3D domains:
\begin{lemma} \label{lemma_E_3d:finite-diff_bound}
   Let $Q\in \mathbb{R}^3$. Then for $\forall \ \epsilon >0$, there exists $\delta>0$ such that 
   $$\displaystyle\sum_{\mathscr{A}(Q_{h}^{+})} h^3 |\sigma_{\alpha x_1 x_2 x_3}|^2 \leq (1+\epsilon) \Big\|\frac{\partial^3 \sigma}{\partial x_1 \partial x_2 \partial x_3} \Big\|^2_{L_2(Q_{h})}$$
whenever $h<\delta$.
\end{lemma}
Although it is more technical, the proof of Lemma~\ref{lemma_E_3d:finite-diff_bound} is very similar to the proof of Lemma~\ref{lemma_E:finite-diff_bound}.
Lemmas~\ref{lemma_E:finite-diff_bound} and ~\ref{lemma_E_3d:finite-diff_bound} imply that the mappings $\mathscr{P}_{h}$ and $\mathscr{Q}_{h}$ satisfy the conditions of Lemma \ref{lem_E:necessary condition}. 
\begin{corollary}\label{prp_E:nec_cond_proof} Assume $Q\in \mathbb{R}^2$ or $\mathbb{R}^3$.
For arbitrary sufficiently small $\epsilon > 0$ there exists $h_{\epsilon}$ such that
\begin{gather}
  \mathscr{Q}_{h}(v) \in \mathscr F^R_{h} \quad \text{
for all}\quad v \in \mathscr F^{R-\epsilon} \quad\text{and}\quad h \leq h_{\epsilon},\label{discrmap}\\ 
\mathscr{P}_{h}([v]_{h}) \in \mathscr F^{R+\epsilon} \quad \text{for all}\quad [v]_{h} \in \mathscr F^R_{h} \quad \text{and}\quad h \leq h_{\epsilon}. \label{interpmap}
\end{gather}
\end{corollary}
To prove \eqref{discrmap}, we first choose $h_{\epsilon}'$ such that for $\forall h<h_{\epsilon}'$ we have
\begin{equation}\label{approxnorm1}
\|\sigma\|_{\tilde{H}^1(Q_h)}^2 \leq \Big (R-\frac{\epsilon}{2}\Big )^2
\end{equation}
Then we apply Lemmas~\ref{lemma_E:finite-diff_bound},~\ref{lemma_E_3d:finite-diff_bound}, Proposition 11 in  \cite{abdulla2019optimal} with $\epsilon_1=\Big (\frac{R}{R-\frac{\epsilon}{2}}\Big )^2-1$, and select $h_\epsilon<h_{\epsilon}'$ such that for $\forall \ h<h_\epsilon$
\begin{equation}
\| \mathscr{Q}_{h}(\sigma)\|_{\mathscr{\tilde H}^1( {Q}_{h})}^2 \leq (1+\epsilon_1)\|\sigma\|_{\tilde{H}^1(Q_h)}^2\leq R^2,
\end{equation}
which proves \eqref{discrmap}. To prove \eqref{interpmap}, we derive the following estimation via straightforward calculation of the respective norm of multilinear interpolation $\sigma'_h$:
\begin{gather}
   \| \mathscr{P}_{h}([\sigma]_{h})\|_{\tilde{H}^1( Q)}^2 = \| \sigma'_{h}\|_{\tilde{H}^1(Q)}^2 
   \leq   \| [\sigma]_{h}\|_{\mathscr{\tilde H}^1( {Q}_{h})}^2 +Ch, \label{approxnorm2}
\end{gather}
where $C$ is independent of $h$. The latter easily imply \eqref{interpmap}.   
Final statement of this section is the following embedding result of \cite{Abdulla-embedding}:
\begin{lemma}\label{embedding} \cite{Abdulla-embedding} 
If $Q\subset \mathbb{R}^2$ or $\mathbb{R}^3$, then 
\begin{equation}\label{embedding1}
 \tilde{H}^1(Q)\hookrightarrow C^{0,\frac{1}{2}}(\overline{Q}); \ \  \ \tilde{H}^1(Q)\Subset C^{0,\mu}(\overline{Q}), 0<\mu < \frac{1}{2}.
 \end{equation}
\end{lemma}
\section{Approximation Theorem and Convergence of the Discrete Optimal Control Problems}\label{section_E:Approximation Theorem}
The following approximation theorem establishes the convergence of the discretized PDE problem:
\begin{theorem} 
\label{thm_E:appro_thm} Let $\{ [v]_{h} \} = \{ ([\sigma]_{h}, U^h) \}$ be a sequence of discrete control vectors such that there exists $R > 0$ for which $[v]_{h} \in \mathscr{F}_{h}^R$ for each $h$, and such that the sequence $\{( \mathscr{P}_{h} ( [\sigma]_{h}), U^h) \}$ converges to $v=(\sigma,U)$ in Tikhonov topology of $\tilde H^1(Q)\times \mathbb{R}^m$ formed with the product of the weak topology of $\tilde{H}^1(Q)$ and Euclidean topology of $\mathbb{R}^m$. Then the sequence of multilinear interpolations $\{u'_{h} \}$ of associated discrete
state vectors $\{[u]_h([v]_h)\}$ converges to the solution $u = u(x;v) \in {H^1}(Q)$ of the elliptic problem \eqref{eq_E:gen_pde}--\eqref{eq_E:robin}, weakly in $H^1(Q)$, strongly in $L_2(Q)$, strongly in $L_2(S)$, and almost everywhere on $Q$. 
\end{theorem}
\textit{Proof.} 
By \eqref{H1energyestimate} of Corollary~\ref{H1boundforu}, sequence $\{u'_h\}$ is uniformly bounded in $ H^1(Q)$. Consequently, $\{u'_{h}\}$ is weakly precompact in $H^1(Q)$. Let $u \in H^1(Q)$ be any weak limit point. By the Rellich-Kondrachev Theorem \cite{nikolskii75}, it is known that there is a subsequence of $\{u'_{h}\}$ that converges to $u$, weakly in $H^1(Q)$, and strongly in $L_{2}(Q)$ and $L_2(S)$. By selecting further subsequence, if necessary, one can achieve that the convergence is almost everywhere on $Q$. We proceed to show that $u$ satisfies the integral identity \eqref{eq_E:weak_sol}.
Without loss of generality, we assume that the whole sequence $\{u'_{h}\}$ converges to $u\in H^1(Q)$. Let $Q'\subset \mathbb{R}^n$ be bounded open domain such that $\bar Q \subset Q'$ and choose arbitrary function $\eta \in C^{1}(\overline Q')$. We assume $h>0$ is small enough that $Q_{h} \subset Q'$. We choose a grid function 
\[ [\eta]_h=\{\eta_\alpha: \ \eta_\alpha=\eta(x_\alpha), \alpha \in \mathscr{A}(Q_h)\} \]
in \eqref{eq_E:disc_weak_sln}. Introducing standard interpolations $\tilde{\eta}_h$ and $\eta_h^i$ as
\[
\tilde{\eta}_{h} |_{C_{h}^{\alpha}} =  \eta_{\alpha}, \ \ \tilde{\eta}_{h}^i |_{C_{h}^{\alpha}} =  \eta_{\alpha x_i}, \quad  \forall \alpha \in \mathscr{A}(Q_{h}^{+})
\]
 we write \eqref{eq_E:disc_weak_sln} in an equivalent form:
\begin{gather}
 \displaystyle\sum_{i=1}^n \int\limits_{Q} \tilde{\sigma}_{h} \tilde{u}_{h}^i \tilde{\eta}_{h}^i\,dx  + \displaystyle\sum_{l=1}^m \frac{1}{Z_l}\int\limits_{E_l}\tilde{u}_{h}   \tilde{\eta}_h\,dS - \displaystyle\sum_{l=1}^m \frac{U^h_l}{Z_l} \int\limits_{E_l}\tilde{\eta}_h\,dS =\nonumber\\
-J_h([u]_h, [\eta]_h)-\displaystyle\sum_{i=1}^n \int\limits_{Q_h\setminus Q} \tilde{\sigma}_{h} \tilde{u}_{h}^i \tilde{\eta}_{h}^i\,dx 
\label{eq_E:disc_weak_sln-1}
\end{gather}
Since, $\tilde{\sigma}_h$ converges to $\sigma$ strongly in $L_2(Q)$, $\tilde{u}^i_h$ converges to $\frac{\partial u}{\partial x_i}$ weakly in $L_2(Q)$, $\tilde{u}_h$ converges to $u$ strongly in $L_2(S)$, $\tilde{\eta}^i_h$ and $\tilde{\eta}_h$ converges to $\frac{\partial \eta}{\partial x_i}$ and $\eta$ uniformly on $\overline{Q}$, the limit of three terms on the left hand side of \eqref{eq_E:disc_weak_sln-1} imply the respective terms of the integral identity  \eqref{eq_E:weak_sol} as $h\to 0$. Hence, it remains to prove that limit of the remaining terms 
in \eqref{eq_E:disc_weak_sln-1} vanishes. By applying Cauchy-Schwartz inequality we have
\begin{equation}\label{Jh1}
|J_h([u]_h, [\eta]_h)| \leq 
 \Big (\displaystyle\sum_{i=1}^n\displaystyle\sum_{\mathscr{A}({Q_{\Delta}^{(i)}})}  h^n u_{{\alpha} x_i}^2\Big )^{\frac{1}{2}}  \| \eta\|_{C^1(\overline{Q'})} \sqrt{nh} \Big ( \displaystyle\sum_{\mathscr{A}(S_h)}  h^{n-1}\Big )^{\frac{1}{2}}. 
  \end{equation}
Noting that every grid point on $S_h$ belongs to cell (with $2^n$ vertices) which intersects $S$, and by recalling the definition of $\hat{S}$ we can estimate 
 \begin{equation}\label{Jh2}
 \displaystyle\sum_{\mathscr{A}(S_h)}  h^{n-1} \leq 2^{n}   \displaystyle\sum_{\mathscr{A}(\hat{S}_h)}h^{n-1}\leq 2^n \sup\limits_{h>0} \displaystyle\sum_{\mathscr{A}(\hat{S}_h)}h^{n-1}=2^n {\cal{H}}^{n-1}(S),
\end{equation} 
where ${\cal H}^{n-1}(\cdot)$ is $n-1$-dimensional Hausdorff measure on $\mathbb{R}^n$.  Since $S$ is Lipschitz, ${\cal{H}}^{n-1}(S)$ coincides with the surface measure $m_{n-1}(S)$ \cite{evansgariepy}. Therefore, from \eqref{Jh1},\eqref{Jh2} and discrete energy estimate \eqref{eq_E:disc_energy_est} it follows that 
\begin{equation}\label{Jh3}
J_h([u]_h, [\eta]_h)=O(\sqrt{h}) \to 0, \ \text{as} \ h\to 0.
\end{equation}
Using Cauchy-Schwartz inequality for the second term in the right hand side of \eqref{eq_E:disc_weak_sln-1} we have
\begin{gather}
\Big |  \displaystyle\sum_{i=1}^n \int\limits_{Q_h\setminus Q} \tilde{\sigma}_{h} \tilde{u}_{h}^i \tilde{\eta}_{h}^i\,dx\Big | \leq \sup\limits_{Q_h\setminus Q}|\tilde{\sigma}_h| \displaystyle\sum_{i=1}^n  \| \tilde u_{h}^i \|_{L_2( Q_{h}\setminus Q)} \| \tilde{\eta}_{h}^i \|_{L_2( Q_{h}\setminus Q)}\nonumber\\
\leq \| \sigma_h'\|_{C(Q_h)} \||[u]_h\||_{\mathscr{ H}^1( {Q}_{h})} \|D \eta\|_{C^1(\overline{Q}')} (m_n(Q_h\setminus Q))^{\frac{1}{2}}\label{qh-q}
\end{gather}
From the embedding result of Lemma~\ref{embedding} and \eqref{approxnorm2} it follows that for sufficiently small $h$
\begin{equation}\label{unifboundsigmah}
\| \sigma_h'\|_{C(Q_h)}\leq C \| \sigma_h'\|_{\tilde{H}^1(Q_h)}\leq C\| [\sigma]_{h}\|_{\mathscr{\tilde H}^1( {Q}_{h})}^2 +1 \leq CR+1.
\end{equation}
Since Lebesgue measure of $Q_h\setminus Q$ converges to zero as $h\to 0$, from the energy estimate \eqref{eq_E:disc_energy_est} and \eqref{unifboundsigmah} it follows that \eqref{qh-q} converges to zero as $h\to 0$. Hence, passing to limit as $h\to 0$, from \eqref{eq_E:disc_weak_sln-1} it follows that the limit function $u$ satisfies the integral identity 
\eqref{eq_E:weak_sol}. \qed

Approximation Theorem~\ref{thm_E:appro_thm} imply the existence of the optimal control.
\begin{corollary}\label{existence}
The optimal control problem  $\mathcal{E}$ has a solution, i.e.
\[
	\mathscr{F}_*:=\Big\{v\in\mathscr{F}^R~\Big|~\mathscr{J}(v)=\mathscr{J}_*\Big\}\neq \emptyset
	\]
\end{corollary}
The proof of the corollary is similar to the proof of existence Theorem 4.4 in \cite{abdulla2018breast}. 

In light of the approximation Theorem~\ref{thm_E:appro_thm},  to complete the proof of Theorem~\ref{thm_E:convergence-of-optimal-control} it remains to prove that the conditions of Lemma~\ref{lem_E:necessary condition} are satisfied. Proof of the condition (iii) of Lemma~\ref{lem_E:necessary condition} coincide with the proof of similar fact from \cite{abdulla13,abdulla2018optimal}. Hence, it only remains to prove that the conditions \eqref{criteriacond1} and \eqref{criteriacond2} of Lemma~\ref{lem_E:necessary condition} are satisfied (see Corollary~\ref{prp_E:nec_cond_proof}).

Let $v \in \mathscr{F}^{(R-\epsilon)}$. By Corollary~\ref{prp_E:nec_cond_proof} we have $ \mathscr{Q}_{h}(\sigma) = [\sigma]_{h} \in \mathscr{F}_{h}^R$. 
Applying Corollary~\ref{prp_E:nec_cond_proof} again, we deduce that $\mathscr{P}_{h} ([\sigma]_{h})$ belong to $ \mathscr{F}^{R+\epsilon}$, and
therefore it forms a weakly precompact sequence in $\tilde{H}^1(Q)$. From compact embedding result of Lemma~\ref{embedding} it follows that it forms a precompact sequence in a strong topology of $C^{0,\mu}(\overline{Q}), \ 0<\mu<\frac{1}{2}$. It easily follows that the whole sequence $\mathscr{P}_{h} ([\sigma]_{h})$ converges to $\sigma$ weakly in $\tilde{H}^1(Q)$, and strongly in $C^{0,\mu}(\overline{Q})$.
From Theorem~\ref{thm_E:appro_thm} it follows that the sequence of multilinear interpolations $\{u'_{h} \}$ of associated discrete
state vectors $\{[u([v]_h)]_h\}$ converges to the solution $u = u(x;v) \in {H^1}(Q)$ of the elliptic problem \eqref{eq_E:gen_pde}--\eqref{eq_E:robin}, weakly in $H^1(Q)$, strongly in $L_2(Q)$ and $L_2(S)$, and almost everywhere on $Q$. Claim (e) of Lemma~\ref{thm_E:interp_relations} implies that the sequence $\tilde{u}_h$ converges to $u$ strongly in $L_2(S)$. Therefore, we have
\begin{gather*}
\lim\limits_{h\to 0} \mathscr{J}_{h}( \mathscr{Q}_{h}(v)=\lim\limits_{h\to 0}\Big ( \displaystyle\sum_{l=1}^m  \Big(  \displaystyle\int_{E_l}  \frac{U_l-\tilde u_{h}}{Z_l} ds  - I_l \Big)^2+\beta |U-U^*|^2 \Big )\\
= \displaystyle\sum_{l=1}^m  \Big(  \displaystyle\int_{E_l}  \frac{U_l-u}{Z_l} ds  - I_l \Big)^2+\beta |U-U^*|^2=\mathscr{J}(v)
\end{gather*}
which proves  \eqref{criteriacond1}.

Let $\{ [v]_{h} =([\sigma]_h,U^h)\}\in  \mathscr{F}^R_{h}$ be arbitrary sequence. From the Corollary~\ref{prp_E:nec_cond_proof} it follows that $(\mathscr{P}_{h} ([\sigma]_{h}),U^h)\in\mathscr{F}^{R+1}$ for sufficiently small $h$, and
therefore it is a precompact sequence in Tikhonov topology of $\tilde{H}^1(Q)\times\mathbb{R}^m$ formed as a product of weak topology of $\tilde H^1(Q)$ and Euclidean topology of $\mathbb{R}^m$. From compact embedding result of Lemma~\ref{embedding} it follows that $\{\mathscr{P}_{h} ([\sigma]_{h})\}$ is a precompact sequence in a strong topology of $C^{0,\mu}(\overline{Q}), \ 0<\mu<1$. Without loss of generality assume that the whole sequence $(\mathscr{P}_{h} ([\sigma]_{h}),U^h)$ converges to some limit $\tilde v=(\tilde \sigma, \tilde U) \in \tilde{H}^1(Q)\times \mathbb{R}^m$. We have
\[ \mathscr{J}( \mathscr{P}_{h}([v]_{h})) - \mathscr{J}_{h}([v]_{h})=\mathscr{J}( \mathscr{P}_{h}([v]_{h})) - \mathscr{J}(\tilde v)+\mathscr{J}(\tilde v) - \mathscr{J}_{h}([v]_{h}) \]
From Theorem~\ref{thm_E:appro_thm} it follows that
\[ \lim\limits_{h\to 0} (\mathscr{J}( \mathscr{P}_{h}([v]_{h})) - \mathscr{J}(\tilde v))=0. \]
The proof of the limit
\[ \lim\limits_{h\to 0} (\mathscr{J}(\tilde v) - \mathscr{J}_{h}([v]_{h}))=0 \]
is almost identical to the preceding proof of \eqref{criteriacond1}. Hence, \eqref{criteriacond2} is proved and this completes the proof of the Theorem~\ref{thm_E:convergence-of-optimal-control}. 

\section{Conclusions}
This paper is on the analysis of the Inverse Electrical Impedance Tomography (EIT) problem on recovering electrical conductivity and potential in the body based on the measurement of the boundary voltages on the $m$ electrodes for a given electrode current. The variational formulation is pursued in the PDE constrained optimal control framework, where electrical conductivity and boundary voltages are control parameters, state vector-potential is a solution of the mixed problem for the second order elliptic PDE, and the cost functional is the norm difference of the boundary electrode current from the given current pattern and boundary electrode voltages from the measurements. The novelty of the control theoretic model is its adaptation to clinical situation when additional "voltage-to-current" measurements can increase the size of the input data from the number of electrodes $m$ up to $m!$ while keeping the size of the unknown parameters fixed. EIT optimal control problem is fully discretized using the method of finite differences. New Sobolev-Hilbert space is introduced, and the convergence of the sequence of finite-dimensional optimal control problems to elliptic coefficient optimal control problem is proved both with respect to functional and control in 2- and 3-dimensional domains. 
\label{sec:conclusions}

}

\begin{thebibliography}{9}

\bibitem{abdulla13}  
    \newblock U.G. Abdulla, 
    \newblock On the optimal control of the free boundary problems for the second order parabolic equations. {I}. {W}ell-posedness and convergence of the method of lines,
    \newblock \emph{Inverse Problems and Imaging}, \textbf{7} (2013), no. 2, 307--340.

\bibitem{abdulla15} 
    \newblock U.G. Abdulla, 
    \newblock On the optimal control of the free boundary problems for the second order parabolic equations. {II}. {C}onvergence of the method of finite differences,
    \newblock \emph{Inverse Problems and Imaging}, \textbf{10} (2016), no. 4, 869--898. 

\bibitem{abdulla2018frechet}  
    \newblock U.G. Abdulla and J.M. Goldfarb,
    \newblock Fr{\'e}chet differentiability in {B}esov spaces in the optimal control of parabolic free boundary problems,
    \newblock \emph{Journal of Inverse and Ill-posed Problems}, \textbf{26} (2018), no. 2, 211--227. 

\bibitem{abdulla2017frechet}  
    \newblock U.G. Abdulla, E. Cosgrove and J. Goldfarb, 
    \newblock On the {Fr\'echet} differentiability in optimal control of coefficients in parabolic free boundary problems,
    \newblock \emph{Evolution Equations and Control Theory}, \textbf{6} (2017), no. 3, 319--344. 
    
\bibitem{abdulla2018optimal}  
    \newblock U.G. Abdulla and B. Poggi, 
    \newblock Optimal control of the multiphase {S}tefan problem,
    \newblock \emph{Applied Mathematics \& Optimization}, \textbf{80} (2019), no. 2, 479-513.  
    
\bibitem{abdulla2019cam}  
    \newblock U.G. Abdulla, V. Bukshtynov and A. Hagverdiyev, 
    \newblock Gradient method in {H}ilbert-{B}esov spaces for the optimal control of parabolic free boundary problems,
    \newblock \emph{Journal of Computational and Applied Mathematics}, \textbf{346} (2019), 84--109.
    
\bibitem{abdulla2020cam}  
    \newblock U.G. Abdulla, J. Goldfarb and A. Hagverdiyev, 
    \newblock Optimal control of coefficients in parabolic free boundary problems modeling laser ablation,
    \newblock \emph{Journal of Computational and Applied Mathematics}, \textbf{372} July 2020, 112736.
    
   
\bibitem{abdulla2019optimal}  
    \newblock U.G. Abdulla,  and B. Poggi,  
    \newblock Optimal {S}tefan problem,
    \newblock \emph{Calculus of Variations and Partial Differential Equations}, \textbf{59}, 61( 2020).     

\bibitem{abdulla2018breast}  
    \newblock U.G. Abdulla, V. Bukshtynov, S. Seif,  
    \newblock Cancer detection through Electrical Impedance Tomography and optimal control theory: theoretical and computational analysis,
    \newblock \emph{arXiv:1809.05936} (2018).
    
    \bibitem{abdulla-cosgrove2020}  
    \newblock U.G. Abdulla and E. Cosgrove,  
    \newblock Optimal control of singular parabolic PDEs modeling {S}tefan-type free boundary problems,
    \newblock \emph{arXiv:2006.07426} (2020).
    
    \bibitem{Abdulla-embedding}
      \newblock U.G. Abdulla,
      \newblock On the Embedding of the Space of Weakly Differentiable Functions into H\"{o}lder Spaces,
      \newblock \emph{preprint} (2020)

\bibitem{alessandrini1988stable}  
    \newblock G. Alessandrini, 
    \newblock Stable determination of conductivity by boundary measurements,
    \newblock  \emph{Applicable Analysis}, \textbf{27} (1988), 153--172.


\bibitem{alsaker}  
    \newblock M. Alsaker and J.L. Mueller, 
    \newblock A $D$-bar algorithm with a priori information for 2-dimensional
electrical impedance tomography,
    \newblock \emph{SIAM J. Imaging Science}, \textbf{9} (2016), 1619--1654.  

 \bibitem{ammari1} 
    \newblock M. Ammari and H. Kang, 
    \newblock Reconstruction of Small Inhomogeneities from Boundary,
    \newblock \emph{Springer}, (2004).    
    
\bibitem{ammari2}  
    \newblock H. Ammari, L. Qiu, F. Santosa and W. Zhang, 
    \newblock Determining anisotropic conductivity using
diffusion tensor imaging data in magneto-acoustic tomography with magnetic induction,
    \newblock \emph{Inverse Problems}, \textbf{34} (2017), 201--224.  

\bibitem{ammari3}  
    \newblock H. Ammari, G.S. Alberti, B. Jin, J.K. Seo and W. Zhang, 
    \newblock The linearized inverse problem in
multifrequency electrical impedance tomography,
    \newblock \emph{SIAM J. Imaging Science}, \textbf{9} (2016), 1525-1551. 

\bibitem{astala1}  
    \newblock K. Astala and L. Palvarinta, 
    \newblock Calderon's inverse conductivity problem in the plane,
    \newblock \emph{Annals of Mathematics}, \textbf{163} (2006), 265--299.
    

    
\bibitem{besov79}  
    \newblock O.V. Besov and V.P. Il'in and S.M. Nikol'skii, 
    \newblock Integral Representations of Functions and Imbedding Theorems,
    \newblock {John Wiley \& Sons}, \textbf{1} (1979).

\bibitem{besov79a}  
    \newblock O.V. Besov and V.P. Il'in and S.M. Nikol'skii, 
    \newblock Integral Representations of Functions and Imbedding Theorems,
    \newblock {John Wiley \& Sons}, \textbf{2} (1979).

\bibitem{calderon} 
    \newblock A.P. Calderon, 
    \newblock On an inverse boundary value problem, in Seminar on Numerical Analysis and Its Applications to Continuum Physics,
    \newblock \emph{Soc. Brasileira de Mathematica, Rio de Janeiro}, (1980), 65--73.

 \bibitem{dodd}  
    \newblock M. Dodd and J. Mueller, 
    \newblock A real-time $D$-bar algorithm for 2d electrical impedance tomography
dat,
    \newblock \emph{Inverse Problems and Imaging}, \textbf{8} (2014), 1013--1031.      

\bibitem{dunlop2015bayesian}   
    \newblock  M. Dunlop and A.M. Stuart,
    \newblock The {B}ayesian formulation of {EIT:} analysis and algorithms,
    \newblock \emph{Inverse Problems and Imaging}, \textbf{10} (2016), 1007--1036.

\bibitem{evansgariepy}   
    \newblock L.C. Evans,  R.F. Gariepy
    \newblock Measure Theory and Finite Properties of Functions,
    \newblock \emph{CRC Press, Taylor \& Francis Group}, (2015).

        

            
     
    

\bibitem{siltanen2}   
    \newblock S.J. Hamilton, M. Lassas and S. Siltanen , 
    \newblock A hybrid segmentation and $D$-bar method
for electrical impedance tomography,
    \newblock \emph{SIAM J. Imaging Science}, \textbf{9} (2016), 770--793.  

\bibitem{siltanen1}   
    \newblock S.J. Hamilton, M. Lassas and S. Siltanen , 
    \newblock A direct reconstruction method for anisotropic
electrical impedance tomography,
    \newblock \emph{Inverse Problems}, \textbf{30} (2014), 770--793.

\bibitem{harrach2019}   
    \newblock B. Harrach, 
    \newblock Uniqueness and Lipschitz stability in electrical impedance tomography with finitely many electrodes,
    \newblock \emph{Inverse problems}, \textbf{19} (2019), 19.  
    
\bibitem{harrah}   
    \newblock B. Harrah and M.N. Minh, 
    \newblock Enhancing residual-based techniques with shape reconstruction
features in Electrical Impedance Tomography,
    \newblock \emph{Inverse Problems}, \textbf{32} (2016).  



\bibitem{holder2004electrical}  
     \newblock D.S. Holder,
     \newblock \emph{Electrical impedance tomography: methods, history and applications},
     \newblock {CRC Press}, 2004.

\bibitem{hyvonen}  
    \newblock N. Hyv\"{o}nen, L. P\"{a}lv\"{a}rinta and J.P. Tamminen,
    \newblock Enhancing D-bar reconstructions for electrical
impedance tomography with conformal maps,
    \newblock \emph{Inverse Problems and Imaging}, \textbf{12} (2018), 373--400.       


\bibitem{jin2}  
    \newblock B. Jin, T. Khan and P. Maass, 
    \newblock A reconstruction algorithm for electrical impedance tomography
based on sparsity regularization,
    \newblock \emph{Int.  J. Numer. Methods}, \textbf{89} (2012), 337--353.  

\bibitem{jin}  
    \newblock B. Jin, Y. Xu and J. Zou, 
    \newblock A convergent adaptive finite element method for electrical impedance
tomography,
    \newblock \emph{IMA J. Numer. Anal.}, \textbf{37} (2017), 1520--1550. 
 
 \bibitem{dbar1}  
    \newblock K. Knudsen, M. Lassas, J.  Mueller and S. Siltanen, 
    \newblock D-Bar method for electrical impedance
tomography with discontinuous conductivities,
    \newblock \emph{SIAM Journal on Applied Mathematics}, \textbf{67} (2007), 893--913.
    
\bibitem{dbar2}  
    \newblock K. Knudsen, M. Lassas, J.  Mueller and S. Siltanen, 
    \newblock Reconstructions of piecewise constant
conductivities by the D-bar method for electrical impedance tomography,
    \newblock \emph{Journal of Physics:
Conference Series}, \textbf{124} (2008).

\bibitem{dbar3}  
    \newblock K. Knudsen, M. Lassas, J.  Mueller and S. Siltanen, 
    \newblock Regularized D-bar method for the
inverse conductivity problem,
    \newblock \emph{Inverse Problems and Imaging}, \textbf{3} (2009), 599--624.

 \bibitem{kaipio2000statistical}  
    \newblock J.P. Kaipio, V. Kolehmainen, E. Somersalo and M. Vauhkonen, 
    \newblock Statistical inversion and {M}onte {C}arlo sampling methods in electrical impedance tomography,
    \newblock \emph{Inverse problems}, \textbf{16} (2000).  

\bibitem{kaipio1999inverse}  
    \newblock J. P. Kaipio, V. Kolehmainen, M. Vauhkonen and E. Somersalo, 
    \newblock Inverse problems with structural prior information,
    \newblock \emph{Inverse problems}, \textbf{15} (1999).

\bibitem{kaipio2005}  
    \newblock J.P. Kaipio and E. Somersalo , 
    \newblock \emph{Statistical and Computational Inverse Problems},
    \newblock {Springer}, 2005.  
     

\bibitem{kenig2007}  
    \newblock C. Kenig, J. Sj\"{o}strand and G. Uhlmann, 
    \newblock The {C}alderon problem with partial data,
    \newblock \emph{Annals of Mathematics}, \textbf{165} (2007), 567--591.

\bibitem{kenig2013}  
    \newblock C. Kenig, J. Sj\"{o}strand and G. Uhlmann, 
    \newblock The {C}alderon problem with partial data on manifolds and aplications,
    \newblock \emph{Analysis and PDE}, \textbf{6} (2013), 2003--2048.

\bibitem{klibanov} 
    \newblock M. V. Klibanov, J. Li and W. Zhang  , 
    \newblock Convexification of electrical impedance tomography with restricted
Dirichlet-to-Neumann map data,
    \newblock \emph{Inverse Problems}, \textbf{35} (2019).

\bibitem{kohn1}  
    \newblock R.V. Kohn, M. Vogelius, 
    \newblock Determining conductivity by boundary measurements,
    \newblock \emph{Comm. Pure Appl. Math.}, \textbf{37} (1984), 289--298.    
        
\bibitem{kohn2}  
    \newblock R.V. Kohn, M. Vogelius, 
    \newblock Determining conductivity by boundary measurements. II. Interior results,
    \newblock \emph{Comm. Pure Appl. Math.}, \textbf{38} (1985), 643--667.     



 \bibitem{dbar4}  
    \newblock V. Kolehmainen, M. Lassas, P. Ola and S. Siltanen, 
    \newblock Recovering boundary shape and conductivity
in electrical impedance tomography, Inverse Problems and Imaging,
    \newblock \emph{Inverse Problems and Imaging}, \textbf{7} (2013), no. 1, 217--242.    


 \bibitem{roininen2014whittle}  
    \newblock S. Lasanen, J.M.L. Huttunen and L. Roininen,
    \newblock Whittle-Matern priors for Bayesian statistical
inversion with applications in electrical impedance tomography,
    \newblock \emph{Inverse Problems and Imaging},(2014), 561-586.  
   
\bibitem{kwon}  
    \newblock O. Kwon, J.K. Seo and J.R. Yoon, 
    \newblock A real-time algorithm for the location search of discontinuous
conductivities with one measurement,
    \newblock \emph{Inverse Problems}, \textbf{} (2002), 201--224.  


\bibitem{ladyzhenskaya2013boundary}  
    \newblock O.A. Ladyzhenskaya,  
    \newblock The boundary value problems of mathematical physics,
    \newblock \emph{Springer Science \& Business Media} (2013).


\bibitem{lassas2009}  
    \newblock M. Lassas, E. Saksman and S. Siltanen,
    \newblock Discretization-invariant Bayesian inversion and Besov
space priors,
    \newblock \emph{Inverse Problems and Imaging}, (2009), 87-122.  

    

    

    
\bibitem{conductivitytumorhigher}  
    \newblock S. Laufer, A. Ivorra and V. Reuter, 
    \newblock Electrical impedane characterization of normal and cancerous human hepatic tissue,
    \newblock \emph{Physiological Measurements}, \textbf{31} (2010), 995--1009.




\bibitem{lechleiter2008}  
    \newblock A. Lechleiter and A. Rieder,
    \newblock Newton regularization for impedance tomography: convergence by local injectivity,
    \newblock \emph{Inverse problems}, \textbf{24} (2008).      
  


\bibitem{bangti}  
    \newblock G. Matthias, J. Bangti and X. Lu, 
    \newblock An analysis of finite element approximation in electrical
impedance tomography,
    \newblock \emph{Inverse Problems}, \textbf{30} (2014).      

    
 \bibitem{nachman1988reconstructions}  
    \newblock A.I. Nachman,
    \newblock Reconstructions from boundary measurements,
    \newblock \emph{Annals of Mathematics}, \textbf{128} (1988), 531--576.   

\bibitem{nachman1995}  
    \newblock A.I. Nachman, 
    \newblock Global uniqueness for a two-dimensional inverse boundary value problem,
    \newblock \emph{Annals of Mathematics}, \textbf{143} (1996), 71--96.
         
 \bibitem{nikolskii75}  
    \newblock S.M. Nikol'skii, 
    \newblock Approximation of Functions of Several Variables and Imbedding Theorems,
    \newblock {Springer-Verlag}, (1975).

%
%
%
%
\bibitem{seo2}  
    \newblock J.K. Seo, J. Lee, S.W. Kim, H. Zribi and E.J. Woo,
    \newblock Frequency-difference electrical impedance
tomography: algorithm development and feasibility study,
    \newblock \emph{Phys. Meas.}, \textbf{29} (2008), 929--941.  
    
\bibitem{seo1}  
    \newblock J.K. Seo and E.J. Woo, 
    \newblock Magnetic resonance electrical impedance tomography,
    \newblock \emph{SIAM Review}, \textbf{53} (2011), 40--68.         

\bibitem{somersalo1992existence}  
    \newblock E. Somersalo, M. Cheney and D. Isaacson, 
    \newblock Existence and uniqueness for electrode models for electric current computed tomography,
    \newblock \emph{SIAM Journal on Applied Mathematics}, \textbf{52} (1992), 1023--1040.    

    
\bibitem{sylvester1987global}  
    \newblock J. Sylvester and, G. Uhlmann, 
    \newblock A global uniqueness theorem for an inverse boundary value problem,
    \newblock \emph{Physiological Measurements}, (1987), 153--169.



\bibitem{vasil1981methods}  
    \newblock F.P. Vasil'ev,  
    \newblock Methods for Solving Extremal Problems. Minimization problems in function spaces, regularization, approximation,
    \newblock \emph{Moscow, Nauka} (1981).
    
\bibitem{widlak}  
    \newblock T. Widlak and O. Scherzer,
    \newblock Hybrid tomography for conductivity imaging,
    \newblock \emph{Inverse Problems}, \textbf{28} (2012).
    

\end{thebibliography}
\end{document}